\documentclass[mathpazo]{csam}
\pdfoutput=1

\usepackage{bm}
\usepackage{xcolor}
\usepackage{cite}
\usepackage[noend]{algpseudocode}
\usepackage[pagewise]{lineno}
\usepackage{algorithmicx,algorithm,setspace,subfigure}
\usepackage{booktabs,multirow}
  
\usepackage{appendix}

\begin{document}
\title{The Double Regularization Method \\ for Capacity Constrained Optimal Transport}


\author[T. Wu et~al.]{Tianhao Wu\affil{1}\comma\affil{2},
      Qihao Cheng\affil{1}, Zihao Wang\affil{3}, Chaorui Zhang\affil{2}, \\ Bo Bai\affil{2}, Zhongyi Huang\affil{1}~and Hao Wu\affil{1}\comma\corrauth}
\address{\affilnum{1}\ Department of Mathematical Sciences,
         Tsinghua University,
         Beijing 100084, P.R. China. \\
         \affilnum{2}\ Theory Lab, Center Research Institute, 2012 Labs, Huawei Technologies Co. Ltd., Hong Kong SAR. \\ 
         \affilnum{3}\  Department of Computer Science and Engineering, Hong Kong University of Science and Technology, Hong Kong SAR. }
\emails{{\tt wuth20@mails.tsinghua.edu.cn} (T.~Wu, shared first authorships), \\ {\tt cqh22@mails.tsinghua.edu.cn} (Q.~Cheng, shared first authorships), {\tt zwanggc@cse.ust.hk} (Z.~Wang), \\ {\tt chaorui.zhang@gmail.com} (C.~Zhang), {\tt ee.bobbai@gmail.com} (B.~Bai), \\ {\tt zhongyih@tsinghua.edu.cn} (Z.~Huang), {\tt hwu@tsinghua.edu.cn} (H.~Wu). }


\begin{abstract}
    Capacity constrained optimal transport is a variant of optimal transport, which adds extra constraints on the set of feasible couplings in the original optimal transport problem to limit the mass transported between each pair of source and sink. Based on this setting, constrained optimal transport has numerous applications, e.g., finance, network flow. However, due to the large number of constraints in this problem, existing algorithms are both time-consuming and space-consuming. In this paper, inspired by entropic regularization for the classical optimal transport problem, we introduce a novel regularization term for capacity constrained optimal transport. The regularized problem naturally satisfies the capacity constraints and consequently makes it possible to analyze the duality. Unlike the matrix-vector multiplication in the alternate iteration scheme for solving classical optimal transport, in our algorithm, each alternate iteration step is to solve several single-variable equations. Fortunately, we find that each of these equations corresponds to a single-variable monotonic function, and we convert solving these equations into finding the unique zero point of each single-variable monotonic function with Newton's method. Extensive numerical experiments demonstrate that our proposed method has a significant advantage in terms of accuracy, efficiency, and memory consumption compared with existing methods. 

\end{abstract}


\maketitle

\section{Introduction}

In this work, we present a novel time-efficient and space-saving method, to solve the discrete Capacity constrained Optimal Transport~(COT) problem~\cite{benamou2015iterative}, which is a discretization of continuous COT~\cite{korman2015optimal} and defined as
\begin{equation}
    \label{eq:COT}
    \begin{gathered}
        \min_{\bm{\gamma}\in \mathbb{R}_+^{N\times M}}\left\langle \bm{C}, \bm{\gamma}\right\rangle, \\
        \text{ s.t. } \bm{\theta} \preceq \bm{\gamma} \preceq \bm{\eta}, \quad  \bm{\gamma}\bm{1} = \bm{u}, \quad  \bm{\gamma}^T\bm{1} = \bm{v}.   
    \end{gathered}
\end{equation}
Here, $\bm{\eta} = \left[\eta_{ij}\right] \in \mathbb{R}_+^{N\times M}$ and $\bm{\theta} = \left[\theta_{ij}\right] \in \mathbb{R}_+^{N\times M}$ are the given capacity upper bound and lower bound matrices, respectively. $\bm{\gamma} = \left[ \gamma_{ij} \right]_{N \times M}$ is the transport plan, whose element $\gamma_{ij}$ denotes the mass transported from position $x_{i}$ to $y_{j}$, and is upper and lower bounded by $\eta_{ij}$ and $\theta_{ij}$ in the COT problem, i.e., $\theta_{ij} \leq \gamma_{ij} \leq \eta_{ij}$. $\bm{A} \preceq \bm{B}$ denotes that every element $a_{ij}$ in matrix $\bm{A}$ is not larger than corresponding $b_{ij}$ in matrix $\bm{B}$. Note that when $\theta_{ij} = 0$ and $\eta_{ij} = 1$ for any $i = 1,\cdots,N$ and $j=1,\cdots,M$, the discretized COT problem~\eqref{eq:COT} is exactly the classical optimal transport problem, in which $\bm{C} = \left[ c_{ij} \right] \in \mathbb{R}_+^{N\times M}$ denotes the cost matrix. \footnote{In this paper, our discussion is general for any $N$ and $M$. For the sake of simplicity, we assume $N = M$ in the rest of the paper. }

Considering the capacity upper and lower bounds in the capacity constrained optimal transport problem, there are a large number of applications in various fields. For network flow~\cite{ford1957primal}, COT can be formulated as a minimum cost maximum flow problem~\cite{wagner1959class,essid2018quadratically,bassetti2020computation}, in which each edge has a capacity. 
For asset pricing and hedging in finance~\cite{tan2013optimal,dolinsky2014martingale,acciaio2016model}, the duality of martingale constrained optimal transport is closely associated with the fundamental theorem of asset pricing. 

Optimal transport theory~\cite{kacher2021comprehensive,hitchcock1941distribution,kantorovich1942translocation} has been successfully applied in different fields~\cite{rubner2000earth,lee2020unbalanced,chen2022computing,chen2018quadratic,engquist2020quadratic,heaton2022wasserstein,pomerleau2015review,shen2021accurate,ye2022optimal}. Consequently, there are numerous algorithms for solving classical optimal transport problem proposed in different perspectives, such as linear programming methods~\cite{yang2021fast}, primal-dual algorithms~\cite{hu2022efficient}, solving Monge-Amp\`{e}re equation~\cite{de2014monge,prins2015least,benamou2000computational,huesmann2019benamou}, proximal block coordinate descent methods~\cite{huxxxxglobal}, reduction and approximation techiniques for high-dimensional distributions~\cite{mengxxxxlarge,mengxxxxsufficient,lixxxxhilbert,zhangxxxxprojection},  Sinkhorn algorithm and its variants~\cite{sinkhorn1967diagonal,cuturi2013sinkhorn,altschuler2017near,alaya2019screening,liao2022fast,liao2022fast2,lin2022efficiency}, etc. However, the large number of constraints added in the COT problem brings out new challenges for computation, and the methods mentioned above can hardly be applied directly to the COT problem. 
A straightforward solution is to formulate the COT problem as a min-cost-max-flow problem on a complete bipartite graph and then solve it with the network flow algorithms, which incurs cubic time complexity according to~\cite{wilf2002algorithms,williamson2019network}. In addition, as is pointed out in~\cite{dong2020study,lin2022efficiency}, solutions obtained by network flow algorithms are indifferentiable. 
There are also several algorithms~\cite{benamou2015iterative,zitouni2007elaboration,chu2020efficient} designed specifically for solving COT problems, among which Iterative Bregman Projection~(IBP) proposed in~\cite{benamou2015iterative} is the most common method. The main idea of IBP is to project the solution onto part of the constraint set alternately based on Dykstra's projection algorithm~\cite{boyle1986method}, which is regarded as a generalization of the Sinkhorn algorithm~\cite{tenetov2018fast}. In this method, the iterative variables are $O(N^2)$ primal variables and the space complexity is inevitably $O(N^2)$. Experimentally, as is shown in Section~\ref{sec:numerical-experiments}, it is impractical to solve capacity constrained optimal transport problems in large-scale scenarios. Therefore, it is necessary to develop a new method for solving COT. 

In this paper, we propose Double Regularization Method~(DRM) to solve the capacity constrained OT problem, inspired by entropic barrier used for the classical OT problem. By introducing two regularization terms, we formulate a double regularized problem to handle both the lower and upper bound constraints in the original COT problem, and approximate the optimal solution of COT. In the double regularized problem, the feasibility of the capacity constraints is naturally guaranteed. Theoretically, as the regularization parameter goes to zero, the optimal solution of the discretized double regularized problem converges to the optimal solution of the discretized COT problem. More importantly, we can solve the duality of the regularized problem with an alternate iteration scheme, which highly reduces both time and space consumption. 
In each alternate iteration of our method, differing from the Sinkhorn algorithm using closed-form update, one of the dual variable vectors is updated by solving several single-variable nonlinear equations based on the objective function of the double regularized problem. A critical observation is that each single-variable nonlinear equation corresponds to a single-variable monotonic function defined on the positive real axis, where the root of each equation exactly equals to the unique zero point of the function. We show that in each iteration, each element in the dual variable vector is updated by finding the unique zero point of a single-variable monotonic function, and can be efficiently handled by Newton's method. Our proposed method needs only $O\left( N \right)$ dual variables in the iterative scheme, while dual variables in existing methods need at least $O(N^2)$ units of memory to the best of our knowledge. Considering that there are two regularization terms in the novel regularized problem, we call this method Double Regularization Method~(DRM). Extensive experiments show that DRM obtains solutions with dozens of times smaller relative error and an order of magnitude faster compared with Iterative Bregman Projection~(IBP). 

The rest of the paper is organized as follows. Section~\ref{sec:DRM} introduces our proposed algorithm, DRM. The numerical results in Section~\ref{sec:numerical-experiments} show the advantages of our algorithm in both accuracy and efficiency. Finally, Section~\ref{sec:conclusion} concludes the work. 

\section{The Double Regularization Method}
\label{sec:DRM}

In this section, we formulate a double regularization problem and propose Double Regularization Method~(DRM) to approximate the optimal solution to the COT problem. For convenience, we assume that the feasible region of the COT problem~\eqref{eq:COT} is non-empty. That is, we assume that
\begin{equation*}
    \bm{\theta}\bm{1} \preceq \bm{u} \preceq \bm{\eta} \bm{1}, \quad \bm{\theta}^T \bm{1} \preceq \bm{v} \preceq \bm{\eta}^T \bm{1}.  
\end{equation*}
By introducing two regularization terms, the double regularized problem writes
\begin{equation}
    \label{eq:DRM-double}
    \begin{aligned}
        \min_{\bm{\gamma} \in \Pi} & \left\langle \bm{C},\bm{\gamma}\right\rangle + \varepsilon \left\langle \bm{\gamma} - \bm{\theta},\ln(\bm{\gamma}-\bm{\theta}) \right\rangle + \varepsilon\left\langle\bm{\eta}-\bm{\gamma},\ln(\bm{\eta}-\bm{\gamma})\right\rangle, \\
        & \Pi = \left\{ \bm{\gamma} \in \mathbb{R}_{+}^{N \times N} \mid \bm{\gamma} \bm{1} = \bm{u}, \bm{\gamma}^{T} \bm{1} = \bm{v}  \right\}, 
    \end{aligned}
\end{equation}
and the regularized problem~\eqref{eq:DRM-double} satisfies the constraint $\theta_{ij} \leq \gamma_{ij} \leq \eta_{ij}$. 

Here, before considering solving the regularized problem~\eqref{eq:DRM-double}, we first show the equivalent conversion between the discretized COT problem~\eqref{eq:COT} and a simplified COT problem in the following lemma. 

\begin{lemma}
\label{lemma:equivalence}
    The discretized COT problem~\eqref{eq:COT} can be equivalently converted to the problem without extra lower bound constraint, which is defined as
    \begin{equation}
        \label{eq:COT-upper}
        \begin{gathered}
            \min_{\bm{\gamma}\in \mathbb{R}_+^{N\times N}}\left\langle \bm{C}, \bm{\gamma}\right\rangle, \\
            \text{ s.t. } \bm{0} \preceq \bm{\gamma} \preceq \bm{\eta}, \quad  \bm{\gamma}\bm{1} = \bm{u}, \quad  \bm{\gamma}^T\bm{1} = \bm{v}. 
        \end{gathered}
    \end{equation}
    Correspondingly, the double regularized problem~\eqref{eq:DRM-double} can be equivalently converted to the following regularized problem
    \begin{equation}
        \label{eq:DRM}
        \begin{aligned}
            \min_{\bm{\gamma} \in \Pi} & \left\langle C, \bm{\gamma} \right\rangle + \varepsilon \left\langle \bm{\gamma}, \ln\bm{\gamma} \right\rangle + \varepsilon \left\langle\bm{\eta}-\bm{\gamma},\ln(\bm{\eta}-\bm{\gamma})\right\rangle, \\
            \Pi & = \left\{ \bm{\gamma} \in \mathbb{R}_{+}^{N \times N} \mid \bm{\gamma} \bm{1} = \bm{u}, \bm{\gamma}^{T} \bm{1} = \bm{v}  \right\},  
        \end{aligned}
    \end{equation}
    which corresponds to the discretized COT problem without lower bound constraint~\eqref{eq:COT-upper}. 
\end{lemma}

The proof of Lemma~\ref{thm:convergence} is shown in the appendix. Thus, for the sake of simplicity, we consider the problem with only an extra upper bound constraint and solve the corresponding regularized problem~\eqref{eq:DRM} in the rest of the paper. Similar to the proof of convergence of entropic regularization for classical OT~\cite{peyre2019computational}, we state the convergence result of the double regularization w.r.t. $\varepsilon$. 

\begin{theorem}
\label{thm:convergence}
    The optimal solution to the double regularized problem~\eqref{eq:DRM} is unique. Let $\bm{\gamma}_{\varepsilon}$ and $\bm{\gamma}^{*}$ be the solutions of \eqref{eq:DRM} and \eqref{eq:COT-upper}, respectively, then $\bm{\gamma}_{\varepsilon}$ converges to $\bm{\gamma}^{*}$ as $\varepsilon$ goes to $0$. 
\end{theorem}

The proof of Theorem~\ref{thm:convergence} is shown in the appendix. An important part in the proof is that $h(x) = x \ln x + ( 1 - x ) \ln (1 - x)$ is a strictly convex function, and thus the regularized problem is still convex. This fact is crucial not only for the proof of the uniqueness of the optimal solution, but also for the design of our algorithm in the following subsection. 

\subsection{Numerical Algorithm}
\label{subsec:numerical-algorithm}

We solve the double regularized problem~\eqref{eq:DRM} to obtain an approximate solution of the capacity constraint OT problem \eqref{eq:COT-upper}. 
Consider the Lagrangian of \eqref{eq:DRM}
\setlength\multlinegap{0pt}
\begin{multline}
    \label{lag}
    L(\bm{\gamma},\bm{\alpha},\bm{\beta})= \sum_{i=1}^{N} \sum_{j=1}^{N} C_{ij} \gamma_{ij} + \varepsilon \sum_{i=1}^{N} \sum_{j=1}^{N} \gamma_{ij} \ln\gamma_{ij} + \varepsilon \sum_{i=1}^{N} \sum_{j=1}^{N} \left( \eta_{ij} - \gamma_{ij} \right) \ln\left( \eta_{ij} - \gamma_{ij} \right) \\
    + \sum_{i=1}^N \alpha_{i} \left(\sum_{j=1}^N \gamma_{ij} - u_i\right) + \sum_{j=1}^N \beta_{j} \left(\sum_{i=1}^N \gamma_{ij}- v_j\right). 
\end{multline}
According to the first-order KKT condition, by taking the derivative of the Lagrangian with respect to each $\bm{\gamma}_{ij}$ and setting $\dfrac{\partial L}{\partial \bm{\gamma}_{ij}}=0$, we obtain the following equation w.r.t. optimal transport plan $\bm{\gamma}$ (for \eqref{eq:DRM}) and corresponding dual variables $\bm{\alpha}, \bm{\beta}$
\begin{equation}
\label{eq:relation}
    \gamma_{ij}=\dfrac{\varphi_i K_{ij} \psi_j \eta_{ij}}{1+\varphi_i K_{ij} \psi_j}
\end{equation}
where $\varphi_i=e^{-\alpha_i/\varepsilon}$,$\psi_j=e^{-\beta_j/\varepsilon}$,$K_{ij}=e^{-C_{ij}/\varepsilon}$. 
Substituting \eqref{eq:relation} into the Lagrangian in \eqref{lag} and setting $\dfrac{\partial L}{\partial \alpha_{i}}=0$ and $\dfrac{\partial L}{\partial \beta_{j}}=0$, we have
\begin{equation*}
\begin{aligned}
    u_{i} &= \sum_{j=1}^N \dfrac{\varphi_i K_{ij} \psi_j \eta_{ij}}{1+\varphi_i K_{ij} \psi_j} = \sum_{j=1}^N \eta_{ij} -  \sum_{j=1}^N \dfrac{\eta_{ij}}{1+\varphi_i K_{ij} \psi_j}, \quad i = 1, 2, \cdots, N, \\
    v_{j} &= \sum_{i=1}^N \dfrac{\varphi_i K_{ij} \psi_j \eta_{ij}}{1+\varphi_i K_{ij} \psi_j} = \sum_{i=1}^N \eta_{ij} - \sum_{i=1}^N\dfrac{\eta_{ij}}{1+\varphi_i K_{ij} \psi_j}, \quad j = 1, 2, \cdots, N,  
\end{aligned}
\end{equation*}
which corresponds to $2N$ single-variable functions 
\begin{equation}
\label{eq:equation-g-phi}
    g_i(\varphi_i)=\sum_{j=1}^N \dfrac{\eta_{ij}}{1+\varphi_i K_{ij} \psi_j} - \sum_{j=1}^N \eta_{ij} + u_i, \quad i = 1,2, \cdots, N, 
\end{equation}
\begin{equation}
\label{eq:equation-f-psi}
        f_j(\psi_j)=\sum_{i=1}^N\dfrac{\eta_{ij}}{1+\varphi_i K_{ij} \psi_j}-\sum_{i=1}^N\eta_{ij}+v_j, \quad j = 1,2, \cdots ,N, 
\end{equation}
whose zero points are dual variables $\bm{\varphi}$ and $\bm{\psi}$ according to \eqref{eq:relation}. Noticing that for $i=1,2,....,N$, 
\begin{equation*}
    g_i'(\varphi_i)=-\sum^{N}_{j=1}\dfrac{ \eta_{ij}\psi_j K_{ij}}{\left(1+\psi_j K_{ij} \varphi_i\right)^2} < 0, \quad f_j'(\psi_j)=-\sum^{N}_{i=1}\dfrac{ \eta_{ij}\varphi_i K_{ij}}{\left(1+\varphi_i K_{ij} \psi_j\right)^2} < 0, 
\end{equation*}
the $2N$ functions are monotonically decreasing. Moreover,
\begin{equation*}
    \lim_{\varphi_{i} \to 0^{+}} g_i(\varphi_i) = u_{i} > 0, \quad \lim_{\varphi_{i} \to +\infty} g_i(\varphi_i) = u_{i} - \sum_{j=1}^N \eta_{ij} < 0. 
\end{equation*}
There are similar results for functions $f_{j}$. Hence, there exists a unique zero point in the interval $(0, +\infty)$ for each single-variable function $g_i$ or $f_j$. The problem of finding the zero points of $g_i(\varphi_i)$ and $f_j(\psi_j)$ can be easily solved by Newton’s method. 

Thus, we propose a novel alternate iteration scheme. The pseudo-code of the proposed method is presented in Algorithm~\ref{alg:2}. 

\begin{algorithm}
\caption{Double Regularization Method for Constraint Optimal Transport}
\label{alg:2}
\begin{algorithmic}[1]
    \Require $\bm{u}$ and $\bm{v}$ of size $N$, $\varepsilon$, $\bm{K}$, maxiter, $\bm{\eta}$, $L$
    \Ensure $\bm{\gamma}$
    \State $\bm{\varphi}^1 \leftarrow \dfrac{1}{N}\bm{1}$, $\bm{\psi}^1 \leftarrow \dfrac{1}{N}\bm{1}$, $m \leftarrow 1$
	\While{$m < maxiter$}
	    \State Solve $g_i\left(\varphi_i\right) = 0$ in  \eqref{eq:equation-g-phi} for $\varphi_i$ with Newton's method. ($i = 1, \cdots, N$) 
	    \State Solve $f_j\left(\psi_j\right) = 0$ in  \eqref{eq:equation-f-psi} for $\psi_j$ with Newton's method. ($j = 1, \cdots, N$) 
	    \State $m \leftarrow m+1$
	\EndWhile
	\For{$i = 1$ to $N$}
		\For{$j = 1$ to $N$}
			\State $\gamma_{ij} \leftarrow \dfrac{\eta_{ij} \varphi^{m}_i K_{ij} \psi^{m}_j}{1+\varphi^{m}_i K_{ij} \psi^{m}_j}$
		\EndFor
	\EndFor
	\State \textbf{return} $\bm{\gamma}$
\end{algorithmic}
\end{algorithm}

In our method, differing from matrix-vector multiplication in the alternate iterations for classical OT, we calculate dual variables by finding the zero points of $g_i$ and $f_{j}$ ($i,j = 1,2,\cdots,N$) with Newton's method in each iteration and update dual variables alternately. The calculation of each $g_i(\varphi_i)$ (or $f_j(\psi_j)$) costs $O(N)$ time since there are $N$ terms in each function in \eqref{eq:equation-g-phi} and \eqref{eq:equation-f-psi}. Thus, the time complexity of DRM is $O\left( N^2 \right)$, which is the same as the Sinkhorn algorithm for classical OT. 

Moreover, there are two vectors $\bm{\varphi}$ and $\bm{\psi}$ used as dual variables, and the number of dual variables is $O\left( N \right)$. While existing methods need at least $O\left( N^2 \right)$ units of memory for dual variables to the best of our knowledge~\cite{benamou2015iterative,chu2020efficient}. Experimentally, our proposed DRM has significant advantages in efficiency compared with IBP. 

\begin{remark}
The log-domain stabilization technique~\cite{chizat2018scaling} can also be aggregated into our DRM algorithm to solve the severe numerical issues caused by small $\varepsilon$. We just add a few steps after line 4 in Algorithm~\ref{alg:2} when the maximum entry of $\bm{\varphi}$ or $\bm{\psi}$ exceeds a given threshold $\tau$:
\begin{equation*}
    \begin{aligned}
        \bm{\alpha} &\leftarrow \bm{\alpha} + \varepsilon \ln\bm{\varphi}, \\
        \bm{\beta} &\leftarrow \bm{\beta} + \varepsilon \ln\bm{\psi}, \\
        \bm{K} &\leftarrow  \mathrm{diag}\left(e^{\bm{\alpha}/\varepsilon}\right)\bm{K} \mathrm{diag}\left(e^{\bm{\beta}/\varepsilon}\right).
    \end{aligned}
\end{equation*}
\end{remark}

\section{Numerical Experiments}
\label{sec:numerical-experiments}


    
    



 
In this section, we conduct numerical experiments to evaluate the performance of DRM algorithm on both 1D case and 2D case. We choose the Wasserstein-2 metric as the objective function for performance comparison. That is, the element of cost matrix $\bm{C}$ is
\begin{equation*}
    C_{ij} = h^2 (i - j)^2 \text{ (for 1D case)}, \quad
    C_{i_1 i_2 j_1 j_2} = h_x^2 (i_1 - j_1)^2 + h_y^2 (i_2 - j_2)^2 \text{ (for 2D case)},
\end{equation*} 
where $h$, $h_x$, and $h_y$ are grid spacing. We compare the performance and computational cost of the constrained capacity optimal transport problem using Gurobi~\cite{optimization2021gurobi}, Iterative Bregman Projection~(IBP)~\cite{benamou2015iterative} and our proposed DRM. 

In all the experiments, unless particularly stated, the regularization parameter $\varepsilon$ is set to $10^{-3}$ for both IBP and DRM. When the relative error of two adjacent iterations $\Delta_{outer}<10^{-5}$, both IBP and DRM are terminated. Specifically, in each outer iteration of DRM, Newton's method is applied for solving~\eqref{eq:equation-g-phi} and \eqref{eq:equation-f-psi}, and is terminated when the relative error of two adjacent Newton inner iterations $\Delta_{inner}<10^{-5}$. 
When Gurobi runs out of memory in large-scale scenarios, we set ``$N/A$'' (short for ``not available'') for the relative error due to lack of ground truth. For IBP, the computational time may exceed 7200s or the memory required may exceed RAM size, which are denoted by ``-''. For each scale of data, the computational time and relative error results are averaged over 50 random experiments. All the experiments are conducted in Matlab R2018b on a platform with 128G RAM, Intel(R) Xeon(R) Gold 5117 CPU @2.00GHz with 14 cores. 

\subsection{1D grid with uniform capacity constraints}

\begin{table}[t]
        \centering
        \renewcommand{\arraystretch}{1.1}
        \tabcolsep 2.5mm
        {
        \begin{tabular}{ccccccc}
        \hline
        \multirow{2}{*}{$\lambda$} & \multirow{2}{*}{$N$}
        & \multicolumn{2}{c}{\texttt{time (sec)}} & \multirow{2}{*}{\texttt{speed-up ratio}} & \multicolumn{2}{c}{\texttt{rel error}}   \\ 
         \cmidrule(r){3-4} \cmidrule(r){6-7}
        & & \textbf{ours} & IBP &  & \textbf{ours} & IBP \\ 
        \hline
        \multirow{4}{*}{5}
        & 1000 & $1.99 \times 10^{0}$ & $3.40 \times 10^{1}$ 
        & $1.71 \times 10^{1}$ 
        & $2.08 \times 10^{-3}$ & $8.75 \times 10^{-3}$ \\
        & 2000 & $7.87 \times 10^{0}$ & $1.90 \times 10^{2}$ 
        & $2.41 \times 10^{1}$ 
        & $2.14 \times 10^{-3}$ & $8.16 \times 10^{-3}$ \\
        & 4000 & $3.76 \times 10^{1}$ & $1.05 \times 10^{3}$ 
        & $2.79 \times 10^{1}$ 
        & $1.98 \times 10^{-3}$ & $6.95 \times 10^{-3}$ \\
        & 8000 & $1.70 \times 10^{2}$ & $4.96 \times 10^{3}$ 
        & $2.92 \times 10^{1}$ & $N/A$ & $N/A$ \\
        \hline
        \multirow{4}{*}{10}
        & 1000 & $4.53 \times 10^{0}$ & $1.28 \times 10^{1}$ 
        & $2.83 \times 10^{0}$ 
        & $3.04 \times 10^{-2}$ & $4.12 \times 10^{-2}$ \\
        & 2000 & $1.65 \times 10^{1}$ & $7.21 \times 10^{1}$ 
        & $4.37 \times 10^{0}$ 
        & $2.84 \times 10^{-2}$ & $4.35 \times 10^{-2}$ \\
        & 4000 & $6.56 \times 10^{1}$ & $3.41 \times 10^{2}$ 
        & $5.20 \times 10^{0}$ 
        & $2.69 \times 10^{-2}$ & $4.53 \times 10^{-2}$ \\
        & 8000 & $3.00 \times 10^{2}$ & $1.59 \times 10^{3}$ 
        & $5.30 \times 10^{0}$  & $N/A$ & $N/A$   \\
        \hline
        \end{tabular}}
        \caption{The 1D random distribution problem with uniform capacity constraints. The comparison between the iterative Bregman projection~(IBP) algorithm and the double regularization method~(DRM) with the different number of grid points $N$. Columns 3-5 are the averaged computational time of the two algorithms and the speed-up ratio of DRM. Column 6-7 are the relative error of the two algorithms compared with the ground truth obtained by Gurobi. }
        \label{table:1D-results-uniform}
\end{table}

First, two discrete distributions are generated on the grid points from the standard uniform distribution of $[0,1]$,
\begin{equation*}
    \bm{u}=\left(u_1,u_2,\cdots,u_N\right), \quad \bm{v}=\left(v_1,v_2,\cdots,v_N\right),
\end{equation*}
and further normalize $\bm{u}$ and $\bm{v}$ such that
\begin{equation*}
    \sum_{i=1}^N u_i=\sum_{j=1}^N v_j=1. 
\end{equation*}
Referring to~\cite{benamou2015iterative}, we use the same upper bound matrix, which is a constant for each entry of the transport plan. Specifically, all the elements in the upper bound matrix are set to $\lambda/N^2$. If $\lambda$ is too small, the entries of the transport plan tend to become the same, or even worse, the feasible set becomes empty. If $\lambda$ is too large (e.g., $\lambda = N^2$), the capacity constraint becomes redundant. Thus, experimentally, $\lambda$ is set to $5$ and $10$ without loss of generality. 

The computational time and relative error of the two algorithms are presented in Table~\ref{table:1D-results-uniform}. We can see that for $\lambda = 5$ and $\lambda = 10$, compared with IBP, our proposed method achieves a smaller relative error with much faster computational speed, respectively. 

\subsection{1D grid with capacity constraints regarding marginals}
\label{subsec:1D-marginal}
In this subsection, referring to~\cite{chu2020efficient}, a different type of upper bound matrix is introduced to guarantee that the feasible set of the COT problem is non-empty, which is
\begin{equation*}
    \bm{\eta}=2\bm{u}\bm{v}^T+\delta\bm{P}. 
\end{equation*}
Here, $2\bm{u}\bm{v}^T$ is the principal term of the upper bound matrix for guaranteeing that the feasible region is non-empty, since $\bm{u}\bm{v}^T$ is in the feasible region. In addition, a random noise term $\bm{P} = \left[P_{ij}\right]_{N\times N}$ is added to demonstrate the generalizability of our proposed algorithm. $\delta$ is an adjustment parameter of the noise matrix. If $\delta$ is too small, the random noise term $\delta \bm{P}$ becomes negligible. If $\delta$ is too large, e.g., $\delta = 1 / \min(P_{ij})$, the capacity constraint becomes redundant. Experimentally, $\delta$ is set to $0.25$, $1$, $4$ in the following experiments. 

\begin{table}[t]
        \centering
        \renewcommand{\arraystretch}{1.1}
        \tabcolsep 2.5mm
        {
        \begin{tabular}{ccccccc}
        \hline
        \multirow{2}{*}{$\delta$} & \multirow{2}{*}{$N$}
        & \multicolumn{2}{c}{\texttt{time (sec)}} & \multirow{2}{*}{\texttt{speed-up ratio}} & \multicolumn{2}{c}{\texttt{rel error}}   \\ 
         \cmidrule(r){3-4} \cmidrule(r){6-7}
        & & \textbf{ours} & IBP &  & \textbf{ours} & IBP \\ 
        \hline
        \multirow{4}{*}{0.25}
        & 1000 & $1.69 \times 10^{0}$ & $5.10 \times 10^{1}$ 
        & $3.02 \times 10^{1}$
        & $1.25 \times 10^{-3}$ & $3.25 \times 10^{-2}$ \\
        & 2000 & $7.10 \times 10^{0}$ & $2.90 \times 10^{2}$ 
        & $4.08 \times 10^{1}$
        & $1.26 \times 10^{-3}$ & $3.24 \times 10^{-2}$ \\
        & 4000 & $3.44 \times 10^{1}$ & $1.56 \times 10^{3}$ 
        & $4.53 \times 10^{1}$
        & $1.30 \times 10^{-3}$ & $3.27 \times 10^{-2}$ \\
        & 8000 & $1.68 \times 10^{2}$ & - 
        & $N/A$ & $N/A$ & $N/A$ \\
        \hline
        \multirow{4}{*}{1.0}
        & 1000 & $1.95 \times 10^{0}$ & $6.69 \times 10^{1}$ 
        & $3.43 \times 10^{1}$
        & $2.28 \times 10^{-3}$ & $2.26 \times 10^{-2}$ \\
        & 2000 & $8.06 \times 10^{0}$ & $2.37 \times 10^{2}$ 
        & $2.94 \times 10^{1}$
        & $2.34 \times 10^{-3}$ & $2.24 \times 10^{-2}$ \\
        & 4000 & $3.60 \times 10^{1}$ & $1.24 \times 10^{3}$ 
        & $3.44 \times 10^{1}$
        & $2.39 \times 10^{-3}$ & $2.25 \times 10^{-2}$ \\
        & 8000 & $1.70 \times 10^{2}$ & $5.62 \times 10^{3}$ 
        & $3.31 \times 10^{1}$ & $N/A$ & $N/A$   \\
        \hline
        \multirow{4}{*}{4.0}
        & 1000 & $3.06 \times 10^{0}$ & $1.73 \times 10^{1}$ 
        & $5.65 \times 10^{0}$
        & $1.49 \times 10^{-2}$ & $2.87 \times 10^{-2}$ \\
        & 2000 & $1.13 \times 10^{1}$ & $9.24 \times 10^{1}$ 
        & $8.18 \times 10^{0}$
        & $1.53 \times 10^{-2}$ & $2.82 \times 10^{-2}$ \\
        & 4000 & $4.55 \times 10^{1}$ & $4.93 \times 10^{2}$ 
        & $1.08 \times 10^{1}$
        & $1.57 \times 10^{-2}$ & $2.79 \times 10^{-2}$ \\
        & 8000 & $2.01 \times 10^{2}$ & $2.30 \times 10^{3}$ 
        & $1.14 \times 10^{1}$ & $N/A$ & $N/A$   \\
        \hline
        \end{tabular}}
        \caption{The 1D random distribution problem with capacity constraints regarding marginals. The comparison between the iterative Bregman projection~(IBP) algorithm and the double regularization method~(DRM) with the different number of grid points $N$. Columns 3-5 are the averaged computational time of the two algorithms and the speed-up ratio of DRM. Column 6-7 are the relative error of the two algorithms compared with the ground truth obtained by Gurobi. }
        \label{table:1D-results}
\end{table}

\begin{figure}[H] 
	\centering 
	\vspace{-0.25cm} 
	\subfigure{
		\includegraphics[width=0.45\linewidth]{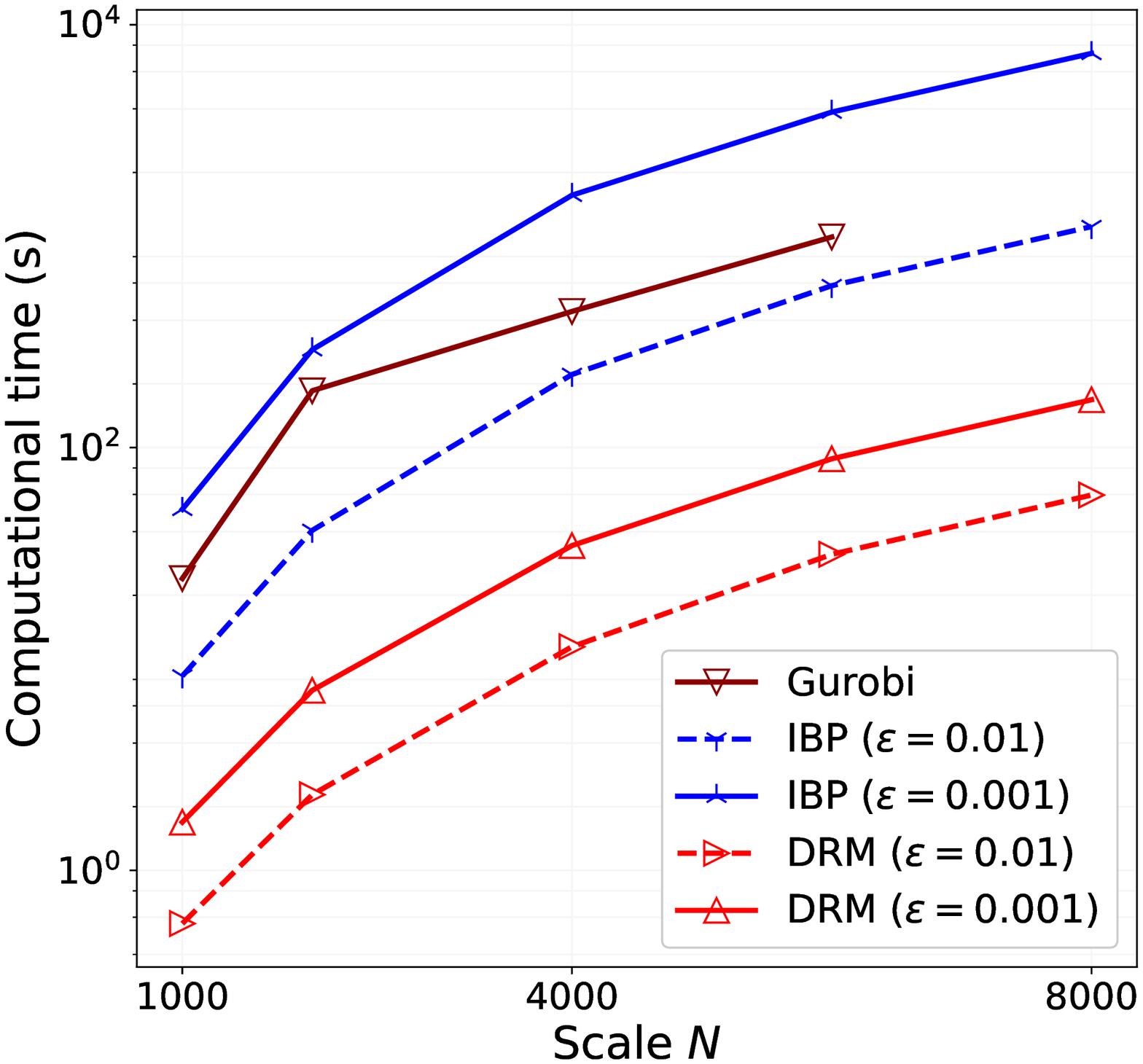}}
	\quad
	\subfigure{
		\includegraphics[width=0.45\linewidth]{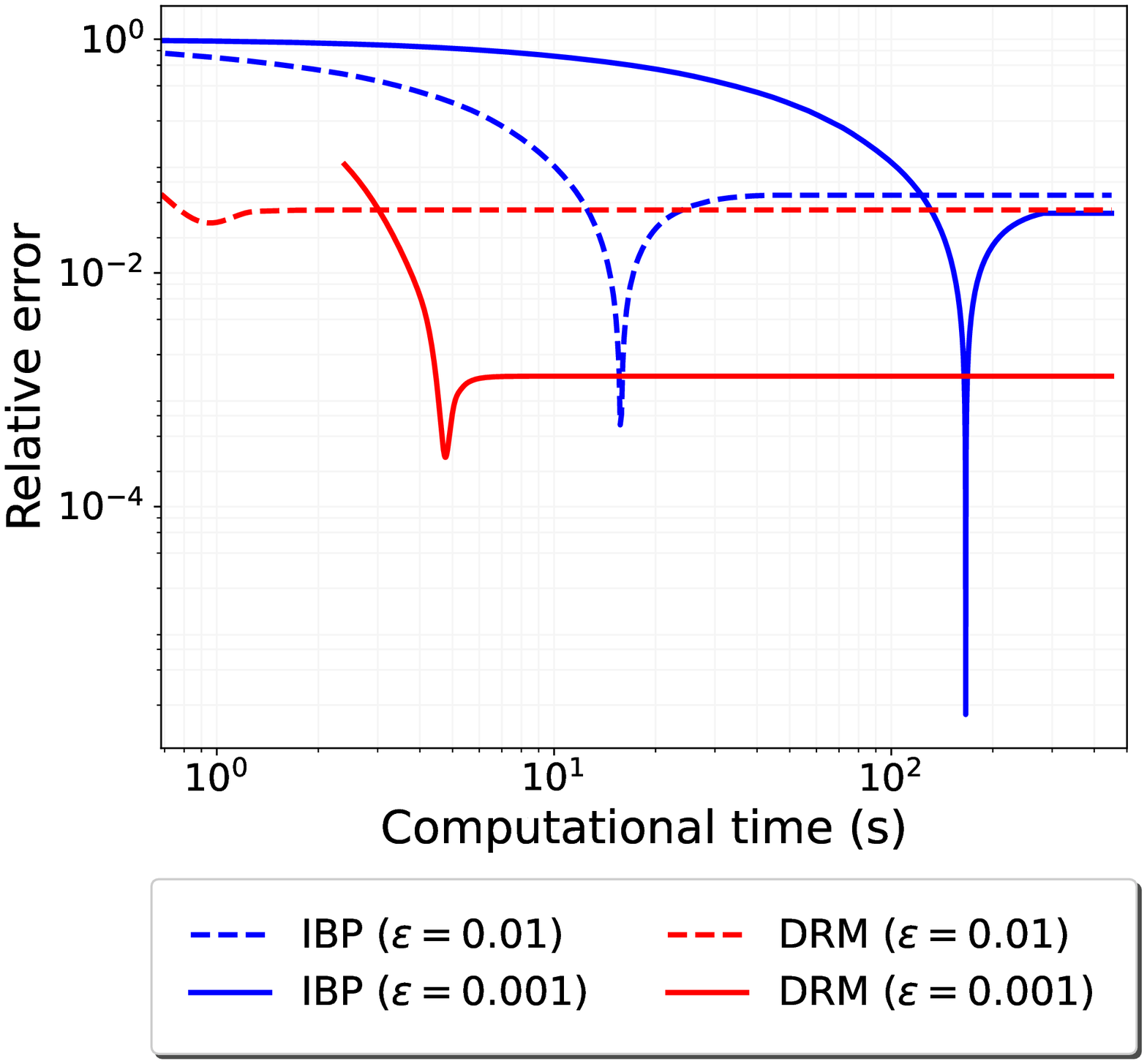}}
	\caption{The 1D random distribution problem with capacity constraints regarding marginals. Left: The comparison of computational time between Gurobi, IBP, and DRM with different numbers of grid points $N$. Right: The relative error between numerical results generated by IBP or DRM and the ground truth w.r.t. computational time. }
	\label{fig:1d-problem}
\end{figure}

The average computational time of DRM and IBP, and the speed-up ratio of DRM on four different sizes of data with $\delta = 0.25$, $1$, and $4$ are presented in Table~\ref{table:1D-results}. From Table~\ref{table:1D-results}, we can see that compared with IBP, DRM has a significant advantage in efficiency for $\delta = 0.25, 1$ and $4$. Particularly, when $\delta$ becomes smaller, the speed-up ratio of DRM becomes larger. This is because when $\delta$ is small, the feasible region is also small, and it is hard for classical algorithms to generate solutions satisfying the feasibility. While our proposed DRM naturally guarantees the feasibility of the capacity constraints. 

Figure~\ref{fig:1d-problem} (left) further demonstrates the algorithm efficiency with the regularization parameter $\varepsilon=0.01$ and $\varepsilon=0.001$ for the three algorithms with different scales of data and $\delta = 0.25$. In Figure~\ref{fig:1d-problem} (right), we discuss the variation of relative error of DRM and IBP w.r.t. computational time when $N=2000$. As is shown in Figure~\ref{fig:1d-problem} (right), for $\delta = 0.25$, the relative error of solutions obtained by DRM with $\varepsilon = 0.001$ is nearly $10^{-3}$, while the relative error obtained by IBP with the same $\varepsilon$ is larger than $10^{-2}$.




\subsection{2D grid with uniform capacity constraints}

Similar to the parameter setting of 1D case with uniform capacity constraints, Wasserstein-2 distance between two generated discrete distributions $\bm{u}=(u_{11},u_{12},....,u_{1N},\cdots,u_{NN})$ and $\bm{v}=(v_{11},v_{12},....,v_{1N},....,v_{NN})$ are computed. All the elements in the two settings of upper bound matrices are $5/N^2$ and $10/N^2$, respectively. 
For the 2D case, the numerical results are shown in Table~\ref{table:2D-results-uniform}. We can see that similar to the 1D case with uniform capacity constraints, our proposed double regularization method~(DRM) also achieves much smaller relative error with dozens of times faster computational speed compared with IBP. 

\begin{table}[t]
        \centering
        \renewcommand{\arraystretch}{1.1}
        \tabcolsep 1.8mm
        {
        \begin{tabular}{ccccccc}
        \hline
        \multirow{2}{*}{$\lambda$} & \multirow{2}{*}{$N \times N$}
        & \multicolumn{2}{c}{\texttt{time (sec)}} & \multirow{2}{*}{\texttt{speed-up ratio}} & \multicolumn{2}{c}{\texttt{rel error}}   \\ 
         \cmidrule(r){3-4} \cmidrule(r){6-7}
        & & \textbf{ours} & IBP &  & \textbf{ours} & IBP \\ 
        \hline
        \multirow{4}{*}{5} 
        & $20 \times 20$   & $2.09 \times 10^{-1}$ & $3.73 \times 10^{0}$ 
        & $1.78 \times 10^{1}$
        & $4.18 \times 10^{-4}$ & $8.71 \times 10^{-3}$   \\
        & $40 \times 40$   & $3.52 \times 10^{0}$  & $9.91 \times 10^{1}$ 
        & $2.82 \times 10^{1}$
        & $4.71 \times 10^{-4}$ & $7.16 \times 10^{-3}$  \\
        & $80 \times 80$   & $1.02 \times 10^{2}$  & $3.49 \times 10^{3}$ 
        & $3.42 \times 10^{1}$ & $N/A$ & $N/A$ \\
        & $160\times160$   & $2.01 \times 10^{3}$  & -  
        & $N/A$  & $N/A$ & $N/A$   \\
        \hline
        \multirow{4}{*}{10} 
        & $20 \times 20$    & $2.12 \times 10^{-1}$ & $2.00 \times 10^{0}$ 
        & $9.43 \times 10^{0}$
        & $2.24 \times 10^{-3}$ & $1.20 \times 10^{-2}$ \\
        & $40 \times 40$    & $3.42 \times 10^{0}$ & $5.18 \times 10^{1}$ 
        & $1.51 \times 10^{1}$
        & $2.61 \times 10^{-3}$ & $8.99 \times 10^{-3}$  \\
        & $80 \times 80$    & $7.06 \times 10^{1}$ & $1.80 \times 10^{3}$ 
        & $2.55 \times 10^{1}$ & $N/A$ & $N/A$ \\
        & $160\times160$    & $1.42 \times 10^{3}$ & - 
        & $N/A$ & $N/A$ & $N/A$  \\
        \hline
        \end{tabular}}
        \caption{The 2D random distribution problem with uniform capacity constraints. The comparison between the iterative Bregman projection~(IBP) algorithm and the double regularization method~(DRM) with the different number of grid points $N \times N$. Columns 3-5 are the averaged computational time of the two algorithms and the speed-up ratio of DRM.  Column 6-7 are the relative error of the two algorithms compared with the ground truth obtained by Gurobi. }
        \label{table:2D-results-uniform}
\end{table}

\begin{table}[t]
        \centering
        \renewcommand{\arraystretch}{1.1}
        \tabcolsep 1.8mm
        {
        \begin{tabular}{ccccccc}
        \hline
        \multirow{2}{*}{$\delta$} & \multirow{2}{*}{$N \times N$}
        & \multicolumn{2}{c}{\texttt{time (sec)}} & \multirow{2}{*}{\texttt{speed-up ratio}} & \multicolumn{2}{c}{\texttt{rel error}}   \\ 
         \cmidrule(r){3-4} \cmidrule(r){6-7}
        & & \textbf{ours} & IBP &  & \textbf{ours} & IBP \\ 
        \hline
        \multirow{4}{*}{0.25} 
        & $20 \times 20$    & $1.70 \times 10^{-1}$ & $4.60 \times 10^{0}$ 
        & $2.71 \times 10^{1}$
        & $1.66 \times 10^{-4}$ & $8.40 \times 10^{-3}$   \\
        & $40 \times 40$    & $2.49 \times 10^{0}$  & $1.18 \times 10^{2}$ 
        & $4.74 \times 10^{1}$
        & $2.08 \times 10^{-4}$ & $8.23 \times 10^{-3}$  \\
        & $80 \times 80$    & $6.12 \times 10^{1}$  & $4.13 \times 10^{3}$ 
        & $6.75 \times 10^{1}$ & $N/A$ & $N/A$ \\
        & $160\times160$    & $1.35 \times 10^{3}$  & -  
        & $N/A$ & $N/A$ & $N/A$  \\
        \hline
        \multirow{4}{*}{1.0} 
        & $20 \times 20$    & $1.62 \times 10^{-1}$ & $3.97 \times 10^{0}$ 
        & $2.45 \times 10^{1}$
        & $2.91 \times 10^{-4}$ & $1.40 \times 10^{-2}$ \\
        & $40 \times 40$    & $2.54 \times 10^{0}$ & $1.01 \times 10^{2}$ 
        & $3.98 \times 10^{1}$
        & $3.44 \times 10^{-4}$ & $1.29 \times 10^{-2}$  \\
        & $80 \times 80$    & $6.18 \times 10^{1}$ & $3.46 \times 10^{3}$ 
        & $5.60 \times 10^{1}$ & $N/A$ & $N/A$ \\
        & $160\times160$    & $1.54 \times 10^{3}$ & - 
        & $N/A$ & $N/A$ & $N/A$   \\
        \hline
        \multirow{4}{*}{4.0} 
        & $20 \times 20$    & $1.79 \times 10^{-1}$ & $2.67 \times 10^{0}$ 
        & $1.49 \times 10^{1}$ 
        & $1.17 \times 10^{-3}$ & $1.59 \times 10^{-2}$   \\
        & $40 \times 40$    & $2.68 \times 10^{0}$  & $6.41 \times 10^{1}$ 
        & $2.39 \times 10^{1}$ 
        & $1.35 \times 10^{-3}$ & $1.39 \times 10^{-2}$ \\
        & $80 \times 80$    & $6.01 \times 10^{1}$  & $2.32 \times 10^{3}$ 
        & $3.86 \times 10^{1}$  & $N/A$ & $N/A$ \\
        & $160\times160$    & $1.37 \times 10^{3}$  & -  
        & $N/A$ & $N/A$ & $N/A$ \\
        \hline
        \end{tabular}}
        \caption{The 2D random distribution problem with capacity constraints regarding marginals. The comparison between the iterative Bregman projection~(IBP) algorithm and the double regularization method~(DRM) with the different number of grid points $N \times N$. Columns 3-5 are the averaged computational time of the two algorithms and the speed-up ratio of DRM. Column 6-7 are the relative error of the two algorithms compared with the ground truth obtained by Gurobi.  }
        \label{table:2D-results}
\end{table}

\subsection{2D grid with capacity constraints regarding marginals}
For the 2D case, we also evaluate the computational cost and performance of DRM, IBP, Gurobi with the same upper bound matrix $\bm{\eta} = 2\bm{u}\bm{v}^T+\delta\bm{P}$ mentioned in Subsection~\ref{subsec:1D-marginal}. 

\begin{figure}[H] 
	\centering 
	\vspace{-0.25cm} 
	\subfigure{
		\includegraphics[width=0.45\linewidth]{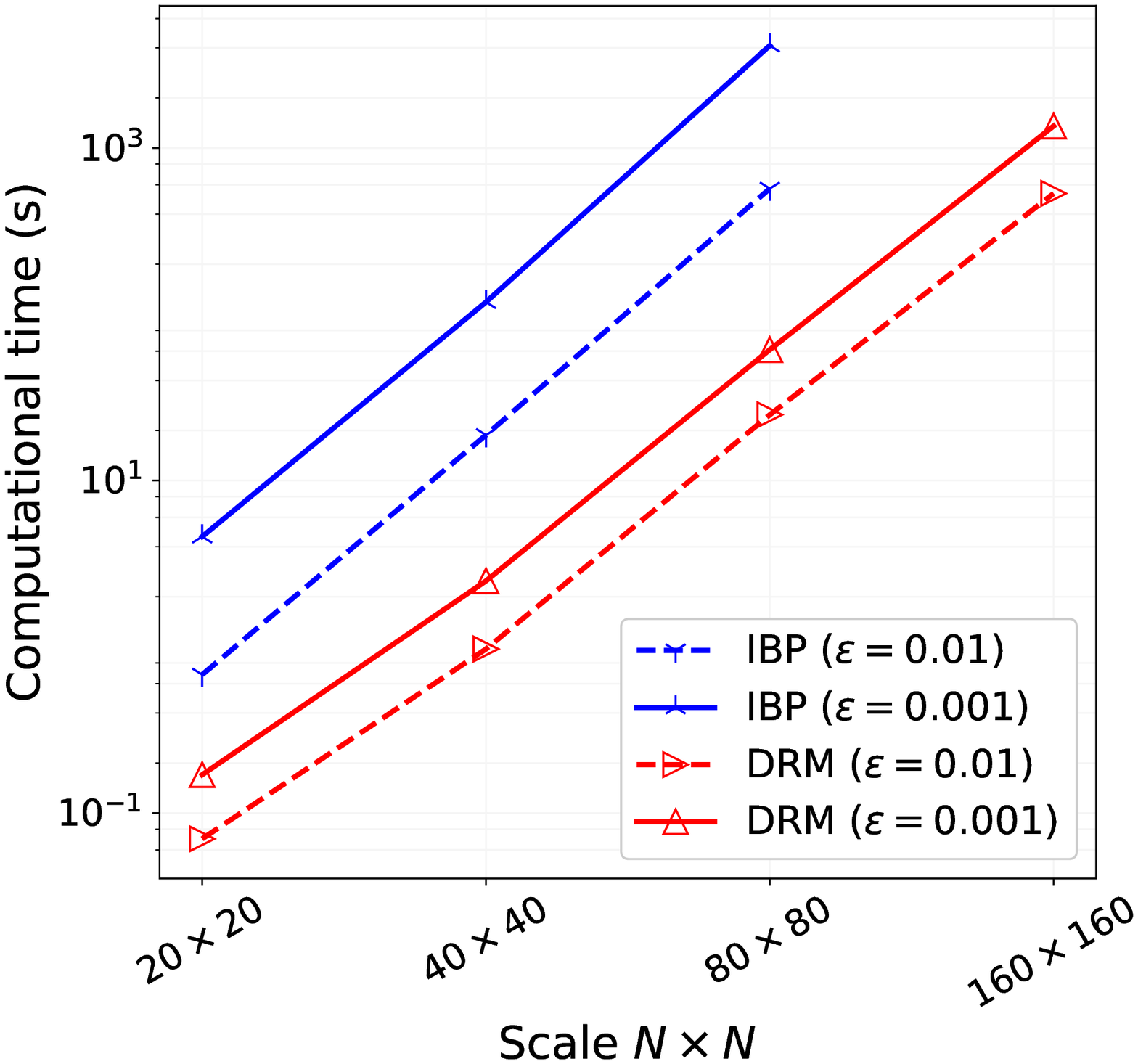}}
	\quad
	\subfigure{
		\includegraphics[width=0.45\linewidth]{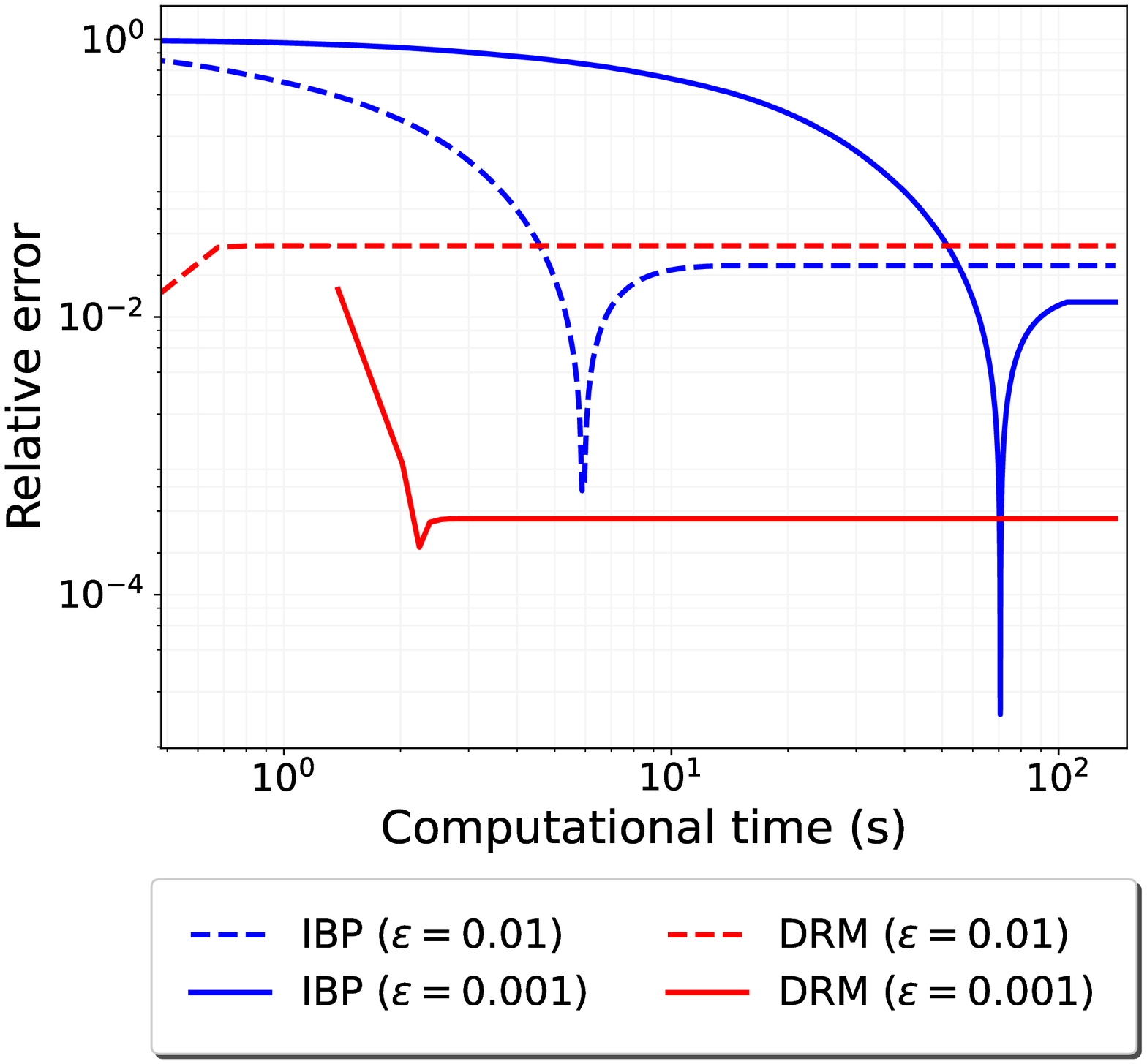}}
	\caption{The 2D random distribution problem with capacity constraints regarding marginals. Left: The comparison of  computational time between IBP and DRM with different number of grid points $N \times N$. Right: The relative error between numerical results generated by IBP or DRM and the ground truth w.r.t. computational time. }
	\label{fig:2d-problem}
\end{figure}

According to Table~\ref{table:2D-results}, we can see that compared with IBP, the computational time of DRM is much smaller. Similar to numerical results in Subsection~\ref{subsec:1D-marginal}, the speed-up ratio becomes larger when $\delta$ becomes smaller. Specifically, with $\delta = 0.25$, DRM is up to $67$ times faster than IBP. Furthermore, Figure~\ref{fig:2d-problem} (left) demonstrates the computational time of the two algorithms with different scales of data. We can see IBP algorithm fails to generate a solution on the $160 \times 160$ data due to the time limitation. At last, Figure~\ref{fig:2d-problem} (right) shows the relative error of solutions obtained by DRM and IBP. Similar to the 1D case, the relative error of DRM with $\varepsilon = 0.001$ and $\delta = 1$ nearly reaches $10^{-4}$, while IBP only generates solutions with relative error $10^{-2}$.  

\subsection{Comparison of memory usage}

\begin{figure}[H] 
	\centering 
	\subfigure{
		\includegraphics[width=0.48\linewidth]{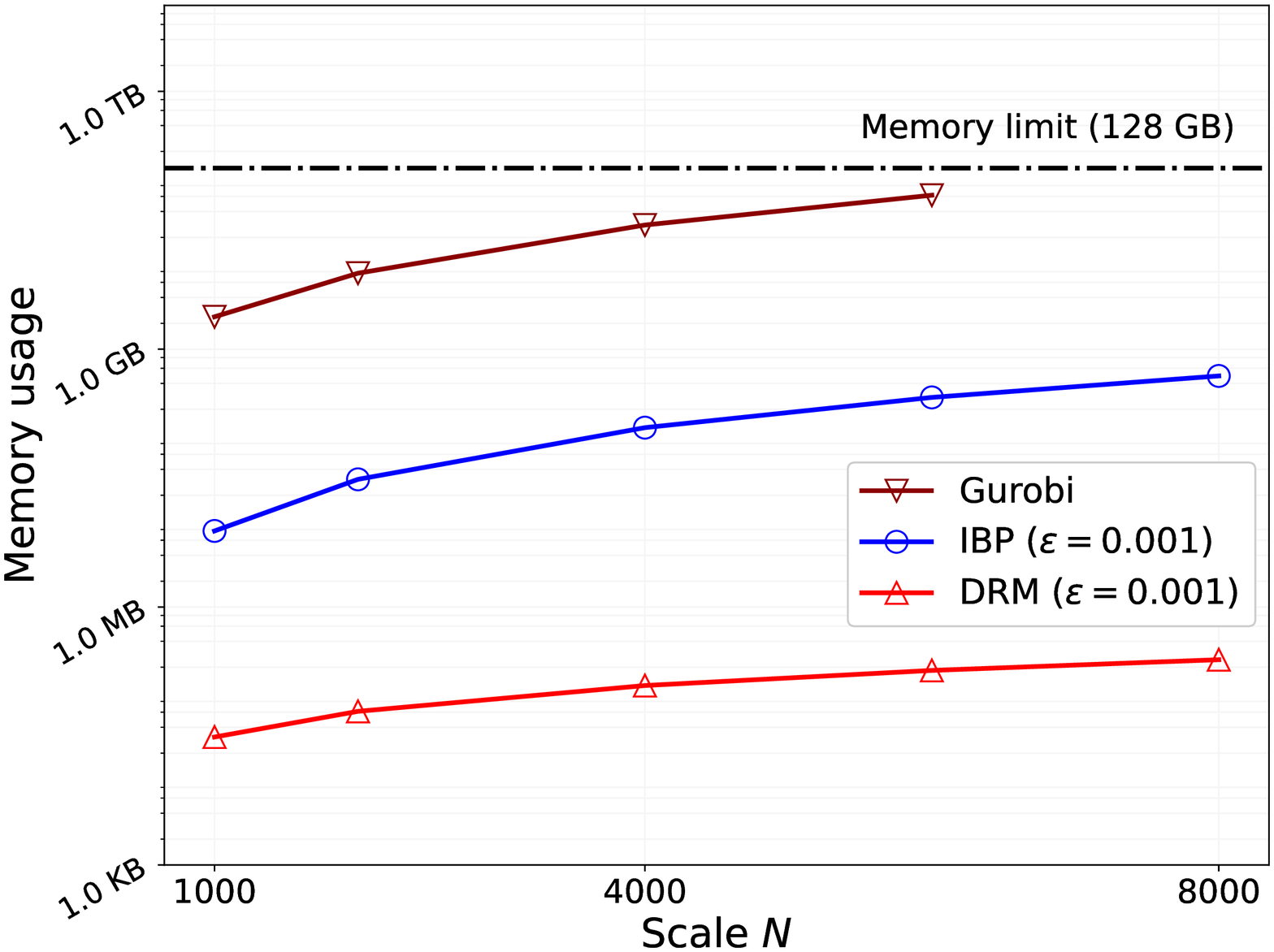}}
	\subfigure{
		\includegraphics[width=0.48\linewidth]{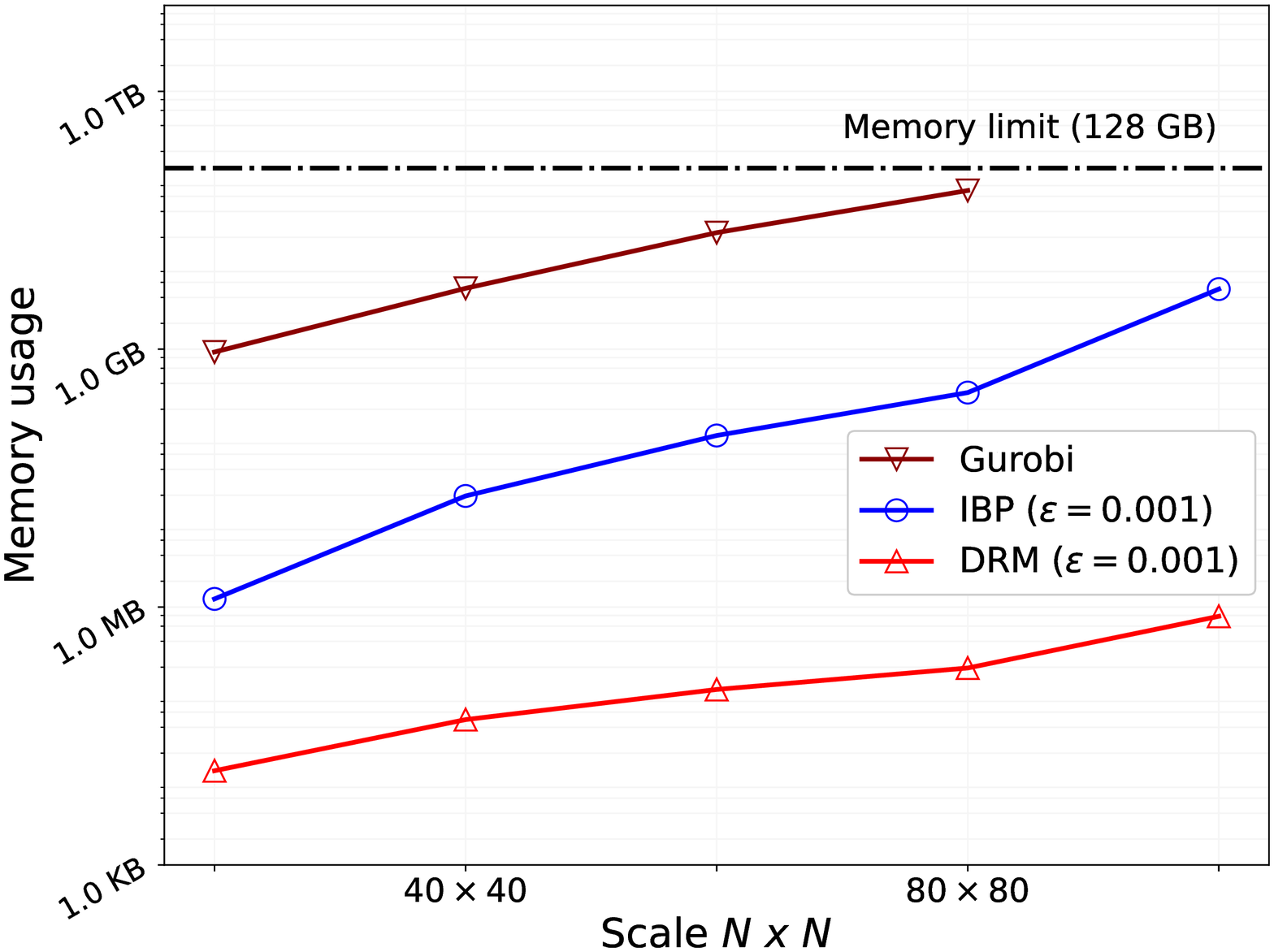}}
	\caption{The memory consumption of Gurobi, IBP and our proposed DRM. Left: The comparison of memory consumption among these three methods with different number of grid points $N$ and upper bound matrix $\bm{\eta} = 2 \bm{u} \bm{v}^T$ in the 1D case. Right: The comparison of memory consumption among these three methods with different number of grid points $N \times N$ and $\bm{\eta} = 2 \bm{u} \bm{v}^T$ in the 2D case. }
	\label{fig:memory-usage}
\end{figure}

Based on the analysis of space complexity of IBP and the proposed DRM in Section~\ref{sec:DRM}, the memory usage of DRM, IBP, and Gurobi is evaluated in this subsection. In the IBP algorithm, the transport plan $\bm{\gamma}$ is modified in each iteration, and thus inevitably incurs $O\left(N^2\right)$ space complexity. On the contrary, our proposed DRM solves the double regularization problem with the alternate iterations of two dual variable vectors, which takes $O(N)$ space. The memory consumption of DRM, IBP, and Gurobi are presented in Figure~\ref{fig:memory-usage}. According to Figure \ref{fig:memory-usage}, DRM takes smaller than $1.0$ MB memory with $N = 8000$ in the 1D case and with $160 \times 160$ in the 2D case. While with the same scales of data, the memory usage of Gurobi is beyond the memory limit (128GB) marked with dash-dotted line in Figure~\ref{fig:memory-usage}. The memory usage of IBP reaches $5$ GB with $160 \times 160$ in the 2D case, and can beyond the memory limit with a larger scale. Compared with the other two methods, the memory consumption of DRM is five orders of magnitude, three orders of magnitude smaller than Gurobi and IBP, respectively. 

\section{Conclusion}
\label{sec:conclusion}

In this paper, we propose Double Regularization Method~(DRM), a novel method for solving the Capacity constrained Optimal Transport problem~(COT). DRM is a generalization of entropy regularization methods for solving classical optimal transport problems. By introducing two regularization terms, we formulate a double regularized problem for COT and make it possible to solve the duality of this problem, which further reduces the memory consumption of variables from $O(N^2)$ to $O(N)$. In the alternate iteration scheme, we update the dual variables by finding zero points of several single-variable monotone functions.  Comprehensive experiments show that our proposed DRM has a remarkable advantage in accuracy, efficiency and memory consumption compared with existing methods. 

\section*{Acknowledgments}

This work was partially supported by National Natural Science Foundation of China Grant
Nos. 11871297, 11871298, 12025104 and 12271289.


\bibliographystyle{abbrv}
\bibliography{ref.bib}

\appendix

\section{PROOFS}
\subsection{Proof of Lemma~\ref{lemma:equivalence}}

Note that the COT problem~\eqref{eq:COT} $\min\limits_{\bm{\gamma} \in \Pi \cap S_{\theta \eta}}\left\langle \bm{C}, \bm{\gamma}\right\rangle$, where $\Pi = \{\gamma_{ij} | \gamma_{ij} \geq 0, \sum_{j=1}^N \gamma_{ij}=u_i, \sum_{i=1}^N \gamma_{ij}=v_j \}$ and $S_{\theta \eta} = \{\gamma_{ij} | \theta_{ij} \leq \gamma_{ij} \leq \eta_{ij} \}$, can be reformulated as
\begin{equation}
\label{eq:reformulate}
    \min_{\bm{\gamma}} \left\langle k_{\theta} \bm{C}, \frac{\bm{\gamma} - \bm{\theta}}{k_{\theta}}\right\rangle
    \end{equation}
    and $\bm{\gamma} = \left[ \gamma_{ij} \right]_{N \times N}$ satisfies
    \begin{equation*}
    \sum_{j=1}^N \frac{\gamma_{ij} - \theta_{ij}}{k_{\theta}} = \frac{ u_i - \sum_{j=1}^{N} \theta_{ij}}{k_{\theta}},\quad \sum_{i=1}^N \frac{\gamma_{ij} - \theta_{ij}}{k_{\theta}} = \frac{v_j - \sum_{i=1}^{N} \theta_{ij}}{k_{\theta}},\quad 0 \leq \frac{\gamma_{ij} - \theta_{ij}}{k_{\theta}} \leq \frac{\eta_{ij} - \theta_{ij}}{k_{\theta}},   
\end{equation*}
where $k_{\theta} = 1 - \bm{1}^{T} \bm{\theta} \bm{1}$ is a positive constant with the fact that $\sum_{i}\sum_{j} \theta_{ij} \leq \sum_{i} u_{i} = 1$. 
Let 
\begin{equation}
\label{eq:transformation}
\begin{aligned}
    & \bm{C}^{\prime} = k_{\theta} \bm{C}, \ \gamma_{ij}^{\prime} = \frac{\gamma_{ij} - \theta_{ij}}{k_{\theta}}, \ \eta_{ij}^{\prime} = \frac{\eta_{ij} - \theta_{ij}}{k_{\theta}}, \\
    & u_{i}^{\prime} = \frac{ u_i - \sum_{j=1}^{N} \theta_{ij}}{k_{\theta}}, \quad v_{j}^{\prime} = \frac{v_j - \sum_{i=1}^{N} \theta_{ij}}{k_{\theta}}. 
\end{aligned}
\end{equation}
Suppose ${\bm{\gamma}}^{*}$ is the optimal solution of the discretized COT problem~\eqref{eq:COT}, then 
\begin{equation*}
    \left\langle \bm{C}, {\bm{\gamma}}^{*} \right\rangle \leq \left\langle \bm{C}, \bm{\gamma} \right\rangle, \quad \forall \bm{\gamma} \in \Pi \cap S_{\theta \eta}, 
\end{equation*}
\begin{equation*}
    \left\langle k_{\theta} \bm{C}, \frac{{\bm{\gamma}}^{*} - \bm{\theta}}{k_{\theta}} \right\rangle \leq \left\langle  k_{\theta} \bm{C}, \frac{\bm{\gamma} - \bm{\theta}}{k_{\theta}} \right\rangle, \quad \forall \bm{\gamma} \in \Pi \cap S_{\theta \eta}. 
\end{equation*}
Considering the substitution in equation~\eqref{eq:transformation}, we have
\begin{equation*}
    \left\langle \bm{C}^{\prime}, \frac{{\bm{\gamma}}^{*} - \bm{\theta}}{k_{\theta}} \right\rangle \leq \left\langle  \bm{C}^{\prime}, \bm{\gamma}^{\prime} \right\rangle, \quad \forall \bm{\gamma}^{\prime} \in \Pi^{\prime} \cap S_{0 \eta^{\prime}}, 
\end{equation*}
and the optimal solution of the discretized COT problem without lower bound constraint~\eqref{eq:COT-upper} is $\dfrac{{\bm{\gamma}}^{*} - \bm{\theta}}{k_{\theta}}$. Hence, \eqref{eq:COT} is equivalently converted to \eqref{eq:COT-upper}.

Then, similar to equation~\eqref{eq:reformulate}, we reformulate the objective function of \eqref{eq:DRM-double}. With the same substitution in equation~\eqref{eq:transformation} and $\varepsilon^{\prime} = k_{\theta} \varepsilon$,  Lemma~\ref{lemma:equivalence} holds. 

\subsection{Proof of Theorem~\ref{thm:convergence}} 

Consider a sequence $\left\{ \varepsilon_{k} \right\}$ where $\varepsilon_{k} > 0$ and $\lim\limits_{k \to \infty} \varepsilon_{k} = 0.$
The optimal solution of the double regularized problem~\eqref{eq:DRM} with $\varepsilon_{k}$ is $\bm{\gamma}_{k}$. Let $S_{0 \eta} = \left\{ \gamma_{ij} | 0 \leq \gamma_{ij} \leq \eta_{ij} \right\}$, then $\bm{\gamma}_{k} \in \Pi \cap S_{0 \eta}$. Considering that $\Pi \cap S_{0 \eta}$ is closed and bounded, there exists a subsequence converging to $\hat{\bm{\gamma}} \in \Pi \cap S_{0 \eta}$. For the sake of simplicity, we denote the subsequence by the same symbol $\left\{ \bm{\gamma}_{k} \right\}$, i.e, $\bm{\gamma}_{k} \rightarrow \hat{\bm{\gamma}}$. For the optimal solution $\bm{\gamma}^{*}$ of \eqref{eq:COT-upper}, with the optimality of $\bm{\gamma}_{k}$ and $\bm{\gamma}^{*}$, we have
\setlength\multlinegap{0pt}
\begin{multline}
    0 \leq \left\langle C, \bm{\gamma}_{k} \right\rangle - \left\langle C, \bm{\gamma}^{*} \right\rangle \leq \varepsilon_{k} \left( \left\langle \bm{\gamma}^{*}, \ln\bm{\gamma}^{*} \right\rangle + \left\langle\bm{\eta}-\bm{\gamma}^{*},\ln(\bm{\eta}-\bm{\gamma}^{*})\right\rangle \right. \\
    \left. - \left\langle \bm{\gamma}_{k}, \ln\bm{\gamma}_{k} \right\rangle - \left\langle\bm{\eta}-\bm{\gamma}_{k},\ln(\bm{\eta}-\bm{\gamma}_{k})\right\rangle \right). 
\end{multline}
Note that the inner product function and logarithmic function are continuous, by taking the limit $k \rightarrow \infty$, we obtain that $\left\langle C, \hat{\bm{\gamma}} \right\rangle = \left\langle C, \bm{\gamma}^{*} \right\rangle$. 

Besides, $h(x) = x \ln x + ( 1 - x ) \ln (1 - x)$ is a strictly convex function, and thus the regularized problem is still convex and has a unique solution converging to $\bm{\gamma}^{*}$ as $\varepsilon$ goes to $0$. 

\end{document}